\documentclass[12pt, a4paper, twoside]{scrartcl}

\usepackage[utf8]{inputenc} 
\usepackage[T1]{fontenc} 
\usepackage{amsmath}
\usepackage{mathtools} 
\usepackage{dsfont}
\usepackage{stmaryrd} 
\usepackage{graphicx} 
\usepackage[english]{babel} 
\RequirePackage[colorlinks,citecolor=blue,urlcolor=blue]{hyperref} 
\usepackage{xcolor}
\usepackage{bm}
\usepackage{amssymb}
\usepackage{enumitem}
\usepackage{wasysym}
\usepackage{fullpage}
\usepackage{soul}

\usepackage{cmbright}
\usepackage[nomath]{lmodern}
\usepackage{inconsolata}
\usepackage{amsthm} 
\theoremstyle{definition}
\newtheorem{thm}{Theorem}[section]
\newtheorem{prop}[thm]{Proposition}
\newtheorem{lem}[thm]{Lemma}
\newtheorem{cor}[thm]{Corollary}
\theoremstyle{remark}
\newtheorem{rem}{Remark}[section]

\usepackage{etoolbox}
\makeatletter
\patchcmd{\@algocf@start}
{-1.5em}
{0pt}
{}{}
\makeatother

\title{The number of connected components in sub-critical random graph processes}
\author{Josué Corujo\footnote{Univ Paris Est Créteil, Univ Gustave Eiffel, CNRS, LAMA UMR 8050, F-94010 Créteil, France}}
\date{ }

\begin{document}

\maketitle

We present a detailed study of the evolution of the number of connected components in sub-critical multiplicative random graph processes.
We consider a model where edges appear independently after an exponential time at rate equal to the product of the sizes of the vertices.
We provide an explicit expression for the fluid limit of the number of connected components normalized by its initial value, when the time is smaller than the inverse of the sum of the square of the initial vertex sizes.
We also identify the diffusion limit of the rescaled fluctuations around the fluid limit.
This is applied to several examples.
In the particular setting of the Erd\H{o}s--Rényi graph process, we explicit the fluid limit of the number of connected components normalized, and the diffusion limit of the scaled fluctuations in the sub-critical regime, where the mean degree is between zero and one.

\section{Introduction}

For $n \in \mathbb{N}$, write $[n]$ for $\{1,\ldots, n\}$.
Let us denote by $G(n,p)$, where $p\in[0,1]$, the Erd\H{o}s-Rényi random graph \cite{Gilbert1959,erdren}: each edge is included in the graph with probability $p$, independently from every other edge.
A continuous-time random graph process related to $G(n,p)$, is naturally constructed as follows: 
fix the vertices set $[n]$ and let each of the $\binom{n}{2}$ edges appear at an exponential time of rate $1$, independently of each other.
This transforms the model into a continuous-time Markov chain, running on the set of graphs with vertices $[n]$ and going from the trivial graph ($n$ disconnected vertices) at time $t = 0$, to the complete graph when $t \to \infty$.
Let us denote this process by $\big( \mathrm{ER}^{(n)}(t) \big)_{t \ge 0}$.
This continuous-time percolation construction is equally obtained by the time-change $t = -\ln(1-p)$ in the natural coupling of $\big(G(n,p), p \in [0,1]\big)$, i.e.\ the process where at each $p \in [0,1]$ the edge $\{i,j\}$ is present in the graph if and only if $\{ U_{\{i,j\}} \le p \}$, for a fixed family $ \displaystyle ( U_{\{i,j\}} )_{i\neq j \in [n]}$ of independent random variables with uniform distribution in $[0,1]$.

A remarkable phenomenon occurs at the scale $t_n = c/n$ where $\big( \mathrm{ER}^{(n)}(t) \big)_{t \ge 0}$ exhibits a phase transition at  $c =1$. 
Namely, with high probability (w.h.p.), for 
$c < 1$, all connected components of the graph are of order 
$O(\ln(n))$, whereas for $c > 1$, there is a single component of order of magnitude $O(n)$ (usually called the \emph{giant component}) and all the others are of order $O(\ln(n))$.
The Erd\H{o}s--Rényi random graph is called sub-critical for $c < 1$ and super-critical for $c > 1$.

The study of the connectivity of the Erd\H{o}s--Rényi random graph $\mathrm{ER}^{(n)}(p)$ has been mostly developed in the regime where $p = \Theta( \ln(n)/n)$.
In this setting it is well understood that if $p = c \ln(n) /n$ with $c < 1$ the graph is disconnected with high probability (w.h.p.), and if $c > 1$ then the graph is connected w.h.p.
However, the connectivity properties of the graph in the sparse regime, i.e.\ $\mathrm{ER}^{(n)}(c/n)$, with $c \ge 0$, has not been so well understood.
In this regime the graph is disconnected w.h.p., so it is interesting to analyze other statistics such as the number of connected components, its limit, and the fluctuations of the convergence toward this limit.
Let $\Rightarrow$ denote  convergence in distribution.
As far as we can recall, the most detailed study of this problem up to this point is carried out by \cite[Thm.\ 2.2]{Puhalskii2005}, where it is proved that for the Erd\H{o}s--Rényi random graph $\mathrm{ER}^{(n)}(c/n)$, we get
\begin{equation}\label{eq:puhalskii}
	\sqrt{n} \left( \frac{1}{n} K_{\mathrm{ER}}^{(n)} \left( \frac{c}{n} \right) - \Big(1 - \frac c 2 \Big) \right)   \Rightarrow \mathrm{N}\left(0 , \frac{c}{2} \right) \text{ for } c < 1 \text{ in } \mathbb{R} \text{ as } n \to \infty,
\end{equation}
where $K_{\mathrm{ER}}^{(n)}(c/n)$ denotes the number of connected components in $\mathrm{ER}^{(n)}(c/n)$,
and $\mathrm{N}(0 , c/2)$ denotes a Gaussian random variable with zero mean and variance $c/2$.
The law of large number result associated to \eqref{eq:puhalskii} is also obtained by Curien \cite[Cor.\ 7.5]{curien2024randomwalkrandomgraphs} by using an exploration process and the differential equation method.
This result is \emph{static} in the sense that it is stated for a fixed value of $c$ but it does not describe the evolution of the trajectories of $\big( \mathrm{ER}^{(n)}(t/n) \big)_{t \ge 0}$.

Enriquez, Faraud and Lemaire \cite{enriquez2023} 
recently considered the super-critical Erd\H{o}s--Rényi random graph process, and studied the evolution of the size of the largest component and the convergence of the sequence of processes of fluctuations.
The same result is re-proved by Lemaire, Limic and the author \cite{CorujoLemaireLimic_2024} by using an exploration process that also allows them to  obtain similar results for the barely-super-critical regime, that occurs before the super-critical regime and after the near-critical regime, to be defined later.
We are interested in obtaining analogous results but for the number of connected components in the sub-critical regime.
Our main result related to the Erd\H{o}s--Rényi random graph process is the following.

Let $D([0,c], \mathbb{R})$ be the space of càdlàg real-valued functions defined on the interval $[0,c]$ and endowed with the Skorohod $J_1$ topology \cite{Skorokhod1956}, see Billingsley \cite{Billingsley1968,Billingsley1999} for background.  

\begin{thm}[Number of CC in the sub-critical Erd\H{o}s--Rényi random graph]\label{thm:ER}
	
	For every $[0, c]$, with $c < 1$, we have that the process
	\[
	\left( \frac{ 1 }{ n } K_{\mathrm{ER}}^{(n)} \Big( \frac{t }{n} \Big) \right)_{t \in [0,c]} \Rightarrow \left( 1 - \frac{t}{2} \right)_{t \in [0,c]} \text{ in } D([0,c], \mathbb{R}) \text{ as } n \to \infty.
	\]
	Furthermore, the process of fluctuations
	\[
	\left( \sqrt{n} \left( \frac{1}{ n }  K_{\mathrm{ER}}^{(n)} \Big( \frac{t }{ n} \Big)  - \Big( 1 - \frac{t}{2} \Big)  \right) \right)_{t \in [0,c]} \Rightarrow \left( B \Big(\frac t 2 \Big) \right)_{t \in [0,c]} \text{ in } D([0,c], \mathbb{R}) \text{ as } n \to \infty,
	\]
	where $B$ is a Brownian motion.
\end{thm}

\begin{rem}[On the mode of convergence]
	In general, the convergence of stochastic processes with càdlàg trajectories is stated in the $J_{1}$ Skorohod topology.
	However, when the limit is continuous, as it is the case in all the convergence results presented in this paper, the convergence holds in the space of càdlàg functions equipped with the uniform topology.
\end{rem}

Notice that this result is consistent with \eqref{eq:puhalskii}, in particular Theorem 2.2 (part 1) in \cite{Puhalskii2005}.
Moreover, the techniques in \cite{Puhalskii2005} rely on the study of an exploration process which does not seem to extend naturally to the study of the dynamical fluctuations.
In contrast, our approach is based on a simple analysis of the martingale problem linked to the evolution of the weighted sizes of the connected components, which readily allows us to establish the convergence of the process of fluctuations.

Rather than the Erd\H{o}s--Rényi random graph process, we study a more general graph process as defined by \cite{Aldous1997} and \cite{Aldous_Limic1998} and recalled in Section \ref{section:model}.
In this model, each vertex has an associated size or mass, and an edge between two vertices appears after an exponential time with rate equal to the product of the sizes. 
We then apply the result to several examples, including the Erd\H{o}s--Rényi random graph.
In particular, Theorem \ref{thm:ER} follows directly from the more general Theorem \ref{thm:main_result}, as well as Corollaries \ref{corol:GER} and \ref{cor:NR}.

\section{The random graph process}
\label{section:model}

Consider $l^2_\searrow$ the space of infinite vectors in $l^2$ whose components are non-negative and ordered in decreasing order and
take $\pmb{\boldsymbol{z}} = (z_1, z_2, \dots) \in l^2_\searrow$ a finite vector, meaning that there exists an integer $\ell \in \mathbb{N}$ such that $z_{\ell} = 0$. 
Aldous \cite{Aldous1997} extended the construction of $\mathrm{ER}^{(n)}$ as follows: instead of mass $1$, let vertex $i\ge 1$ have initial mass $z_i>0$. 
For each $i,j\in [n]$ let the edge between $i$ and $j$ appears after an exponential time with rate $z_i \cdot z_j$, independently of others.
In the sequel, we 
denote by ${G}^{\boldsymbol{z}} = ({G}^{\boldsymbol{z}}(t), t\ge 0)$ a graph-valued continuous-time Markov process following these dynamics.
Note that $\mathrm{ER}^{(n)}$ corresponds to the case where $\boldsymbol{z}$ is the vector with $n$ components equal to one and the other components equal to zero.

Let us denote by $X^{\boldsymbol{z}}(t) \in l^2_\searrow$ the vector recording the weighted sizes of the connected components in $G^{\boldsymbol{z}}(t)$, for $t \ge 0$.
In particular, $X^{\boldsymbol{z}}(0) = \boldsymbol{z}$.
Due to elementary properties of independent exponential random variables, it is immediate that a pair of connected components merges at a rate equal to the product of their masses.
In other words, the vector of weighted sizes of the connected components of the  continuous-time random graph evolves according to the \emph{multiplicative coalescent} (MC) dynamics:
\begin{equation}
	\label{merge}
	\begin{array}{c}
		\mbox{ any pair of components with masses (sizes) $x$ and $y$ merges }\\
		\mbox{ at rate $x \cdot y$ into a single component of mass $x+y$.}
	\end{array}
\end{equation}
Notice that the sizes of the vertices and the dynamics are related as follows: for every $\alpha > 0$ we get
\begin{equation}\label{eq:MC_change_time}
	\Big(X^{ \alpha \boldsymbol{z}} (t), t \ge 0 \Big) \stackrel{\mathcal{L}}{=} \left( \alpha X^{\boldsymbol{z}} \left(\alpha^2 \, t \right), t \ge 0 \right),
\end{equation}
where $\alpha \boldsymbol{z}$ denotes the usual product of a constant by a vector, and $\stackrel{\mathcal{L}}{=}$ denotes the equality in law.

The generator $\mathcal{Q}$ of the MC acts on a test function $g: l^2_\searrow \to \mathbb{R}$ as follows
\[
\big( \mathcal{Q} g \big) (\pmb{x}) = \sum_{i < j} x_i x_j \Big( g(\pmb{x}^{i,j}) - g(\pmb{x}) \Big),
\]
where $\pmb{x} \in l^2_{\searrow}$ and $\pmb{x}^{i,j}$ is the configuration obtained from $\pmb{x}$ by merging the $i^{\mathrm{th}}$ and $j^{\mathrm{th}}$ clusters, and instantaneously reordering the components in decreasing order.

Because this natural relation between $({G}^{\boldsymbol{z}}(t), t\ge 0)$ and the MC dynamic, we will call this random graph process the \emph{multiplicative random graph process}.

Let us assume $\boldsymbol{z}$ is finite, meaning $\kappa(\boldsymbol{z}) < \infty$, where
\[
\kappa : \pmb{x} \in l^2_\searrow \mapsto \sum_{i \ge 1} \mathbf{1}_{(0, \infty)}(x_i).
\]
Applying the generator of the MC, denoted $\mathcal{Q}$ to $g(X^{\boldsymbol{z}}(t))$, where $g$ is an arbitrary function from $l^2_\searrow$ to $\mathbb R$, one can conclude that the process $(M_g(t))_{t \ge 0}$ where
\[
M_g(t) := g\big( X^{\boldsymbol{z}}(t) \big) - g(z) - \int_0^t \big( \mathcal{Q} g \big) \left( X^{\boldsymbol{z}}(s) \right) \mathrm{d}s
\]
is a local $(\mathcal{F}_t)$-martingale, where $\mathcal{F}_t := \sigma (X^{\boldsymbol{z}}(s), s \le t)$.
Moreover, the predictable quadratic variation of $M_g$ satisfies
\[
\langle M_g \rangle_t = \int_0^t \Big( \Gamma g \Big) \big( X^{\boldsymbol{z}}(s) \big) \mathrm{d}s,
\]
where $\Gamma$ is the \textit{carré-du-champ} operator defined as
\begin{equation}\label{eq:def_carre-du-champ}
	\Gamma g := \mathcal{Q} g^2 - 2 g \mathcal{Q} g,
\end{equation}
see for instance \cite[\S\ 15.5.1]{DelMoral-Penev2017}.
Since we are interested in the number of connected components,
let us define
\(
K^{\boldsymbol{z}}(t) := \kappa\big( X^{\boldsymbol{z}}(t) \big),
\)
which is precisely the number of connected components in $G^{\boldsymbol{z}}(t)$.
Define $\sigma^{(k)}(\pmb{x}) := \sum_{i \ge 1} x_i^k$, for $k = 1,2$.
Note that
\begin{equation}\label{eq:calcul_carre-du-champ_kappa}
	\big( \mathcal{Q} \kappa \big) (\pmb{x}) = \sum_{i < j} x_i x_j \Big( \kappa(\pmb{x}^{i,j}) - \kappa(\pmb{x}) \Big) = - \sum_{i < j} x_i x_j = -\frac{1}{2} \Big( \big( \sigma^{(1)} (\pmb{x}) \big)^2 - \sigma^{(2)} (\pmb{x}) \Big).
\end{equation}
The coalescent dynamics preserve the total mass, so
\[
\sum_{i \ge 1} X^{\boldsymbol{z}}(t) = \sigma^{(1)}(\boldsymbol{z}), \text{ for every } t \ge 0.
\]
Furthermore,
\begin{align}
	\Big( \mathcal{Q} \kappa^2 \big) (\pmb{x}) &= \sum_{i < j} x_i x_j \Big( \kappa(\pmb{x}^{i,j})^2 - \kappa(\pmb{x})^2 \Big) \nonumber \\
	&= \sum_{i < j} x_i x_j \Big( \big( \kappa(\pmb{x}) - 1 \big)^2 - \kappa(\pmb{x})^2 \Big) \nonumber \\
	&= - (2 \kappa(\pmb{x}) - 1 ) \sum_{i < j} x_i x_j \nonumber \\
	&= (2 \kappa(\pmb{x}) - 1 ) \Big( \mathcal{Q} \kappa \big) (\pmb{x}). \label{eq:calcul_carre-du-champ_kappa^2}
\end{align}
Hence, plugging  \eqref{eq:calcul_carre-du-champ_kappa^2} into \eqref{eq:def_carre-du-champ} we get
\[
\Gamma \kappa  =  (2 \kappa - 1) \cdot  \mathcal{Q} \kappa - 2 \kappa \cdot  \mathcal{Q} \kappa = - \mathcal{Q} \kappa.
\]
We thus get the next result as an immediate consequence of the previous discussion.

\begin{prop}[The associated martingale problem]\label{prop:martingale_problem}
	The process $M^{\boldsymbol{z}} = \left( M^{\boldsymbol{z}} \left( t \right) \right)_{t \ge 0}$, where
	\[
	M^{\boldsymbol{z}} \left( t \right) = K^{\boldsymbol{z}} \left( t \right) - \kappa(\boldsymbol{z}) + \frac{t}{2}  \big( \sigma^{(1)}(\boldsymbol{z}) \big)^2 - \frac{1}{2} \int_0^t \sum_{i \ge 1} \big( X_i^{\boldsymbol{z}} \big)^2 \left( s \right) \mathrm{d}s,
	\]
	is a local $(\mathcal{F}_t)$-martingale, with predictable quadratic variation
	\[
	\langle M^{\boldsymbol{z}} \rangle_t = \frac{t}{2} \big( \sigma^{(1)}(\boldsymbol{z}) \big)^2 - \frac{1}{2} \int_0^t \sum_{i \ge 1} \big( X_i^{\boldsymbol{z}} \big)^2 \left( s \right) \mathrm{d}s.
	\]
\end{prop}

\begin{rem}[On the quadratic variation]
	Notice that from the expressions of $M^{\boldsymbol{z}}$ and its predictable quadratic variation, we can easily obtain the quadratic variation of $M^{\boldsymbol{z}}$, which satisfies
	\[
	\left[ M^{\boldsymbol{z}} \right]_t =  \kappa(\boldsymbol{z}) - K^{\boldsymbol{z}} \left( t \right).
	\]
	This can also be deduced from the well-known formula for the quadratic variation $\left[ M^{\boldsymbol{z}} \right]_t  = \sum_{0 \le s \le t} (\Delta M_s) ^2$, where $\Delta M_s = M_s - M_{s-}$.
\end{rem}

A classical result in the multiplicative coalescent (MC) literature establishes the following upper bound on the second moment of the MC, see \cite[\S 2.1]{VLthesis} and \cite[Lemma 3.1]{Limic-Vitalii-moments-MC}, for instance.

\begin{lem}[Bound on the second moment of the MC]\label{lemma:bound_2nd-moment_MC}
	For every $t \in \big[ 0, 1/\sigma^{(2)}(\boldsymbol{z}) \big)$,
	\[
	\mathbb{E} \left[ \sum_{i \ge 1} X_i^{\boldsymbol{z}}(t)^2 \right] \le \frac{\sigma^{(2)}(\boldsymbol{z})}{ 1 - t \, \sigma^{(2)}(\boldsymbol{z})}.
	\]
\end{lem}
Although this simple bound on the expectation in Lemma \ref{lemma:bound_2nd-moment_MC} fails at $t = 1/\sigma^{(2)}(\boldsymbol{z})$, the expectation stays finite for every $t \ge 0$, cf.\ \cite{Limic-Vitalii-moments-MC}.

Lemma \ref{lemma:bound_2nd-moment_MC} allows us to easily control the integral terms in the expressions of the martingale \( M^{\boldsymbol{z}} \) and its quadratic variation. This is one of the main reasons we work within this regime.  

Previous studies on the size of the connected components in the multiplicative random graph process \cite{Aldous1997,Aldous_Limic1998,Limic2019} consider what is known as the near-critical regime, defined as  
\[
q_n(t) = \frac{1}{\sigma^{(2)}(\boldsymbol{z}^{(n)})} + t, \quad t \in \mathbb{R},
\]
where \( \boldsymbol{z}^{(n)} \) is a sequence of initial vector sizes satisfying certain conditions (see Section 5 in \cite{Limic2019}).  
It has been shown that in this regime, $q_n(t)$, for every $t \in \mathbb{R}$, the weighted size of the $k^{\mathrm{th}}$ largest connected component of $G^{\boldsymbol{z}^{(n)}}$ is \( \Theta(1) \) as \( n \to \infty \), for every fixed $k \ge 1$, cf.\ \cite[Prop.\ 4]{Aldous1997} and \cite[Cor.\ 10]{Limic2019}. 
Additionally, these sizes become \( o(1) \) as \( t \to -\infty \) and the size of the largest connected component goes to infinity as \( t \to +\infty \).  
In other words, the process recording the weighted sizes of the connected components starts from dust as \( t \to -\infty \) and transitions to a state where a giant component first appears as \( t \to +\infty \). 
Consequently, it is natural to define the earlier times, i.e.,  
\(
{c}/{\sigma^{(2)}(\boldsymbol{z}^{(n)})}, c \in [0,1),
\)
as corresponding to the sub-critical regime.  

Specifically, for the Erd\H{o}s--Rényi random graph taking \( \boldsymbol{z}^{(n)} \) as the vector whose first \( n \) components are equal to \( 1/n^{2/3} \) and zero otherwise (again, the above-mentioned conditions \cite[Section 5]{Limic2019} are important for this choice), and using \eqref{eq:MC_change_time} we get
$$
\left( X^{n^{2/3} \boldsymbol{z}^{(n)}} \left( t \right), t \ge 0 \right) \overset{\mathcal{L}}{=} \left( n^{2/3} X^{\boldsymbol{z}^{(n)}} ( n^{4/3} \, t ), t \ge 0 \right).
$$
The left-hand side of the previous display is the process recording the ordered component sizes in the Erd\H{o}s--Rényi random graph process.
Hence, we recover the well-known result that in the near-critical regime  
$p_n(t) = {1}/{n} + t/n^{4/3}$,
the connected components of $\mathrm{ER}^{(n)}(p_n(t))$, have sizes of order $n^{2/3}$, while in the sub-critical regime  
\(
{c}/{n}, c \in [0,1),
\)
its connected components are of size $o(n^{2/3})$. 
Of course, more precise results on the size of the connected components in the sub-critical regime are well known; see, for instance, \cite{bol-book,Janson2008,Probabilistic_method,Remco2017}.

\subsection{Main result}

Let us consider a sequence of vectors of initial masses $\boldsymbol{z}^{(n)}$ and, for every time $t \ge 0$, its respective random graph $G^{(n)}(t)$, the vector of sizes of its connected components $X^{(n)}(t)$ and the number of connected components $K^{(n)}(t)$.
Furthermore, to avoid notation clutter, let us also denote $\kappa_{n} = \kappa(\boldsymbol{z}^{(n)})$ and $\sigma_n^{(k)} = \sigma^{(k)}(\boldsymbol{z}^{(n)})$, for $k = 1,2$.
We are now able to establish our main result.

\begin{thm}[Limit behavior of the number of CC] \label{thm:main_result}
	Let us assume that $\kappa_{n} < \infty$, for every $n \in \mathbb{N}$, and $\kappa_{n} \to \infty$, and
	there exist two sequences $(\varkappa_{n})_{n}$ and $(\varsigma_{n})_n$ such that
	\[
	\frac{ \kappa_{n} }{ \varkappa_{n} } \xrightarrow[n \to \infty]{} 1 \text{ and } \frac{\sigma_n^{(2)}}{\varsigma_{n}} \xrightarrow[n \to \infty]{} 1.
	\]
	In addition, suppose there exists a non-negative constant $\alpha$ such that
	\begin{equation}\label{eq:cond-Thm}
		\frac{ (\sigma_n^{(1)})^2}{\varkappa_{n} \, \varsigma_{n} }  \xrightarrow[n \to \infty] {} \alpha.
	\end{equation}
	Then, for every $[0,c] \subset [0,1[$ when $n \to \infty$, the sequence of  processes
	\[
	\left( \frac{ 1 }{ \varkappa_{n} } K^{(n)} \left( \frac{t }{ \varsigma_{n} } \right)  \right)_{t \in [0,c]} \Rightarrow \Big( 1 - t \frac{ \alpha}{2} \Big)_{t \in [0,c]} \text{ in } D([0,c], \mathbb{R}) \text{ as } n \to \infty.
	\]
	
	In addition, if there exists two non-negative real constants $\beta_1$ and $\beta_2$ such that
	\begin{equation}\label{eq:cond-Thm2}
		\sqrt{\varkappa_{n}} \left( \frac{\kappa_{n}}{\varkappa_{n}} - 1 \right) \xrightarrow[n \to \infty]{} \beta_1 \; \text{ and } \;
		\sqrt{\varkappa_{n}} \left( \frac{ (\sigma_n^{(1)})^2}{\varkappa_{n} \varsigma_{n}} - \alpha \right) \xrightarrow[n \to \infty] {} \beta_2. 
	\end{equation}
	then, when $n \to \infty$, the sequence of processes
	\[
	\left( \sqrt{\varkappa_{n}} \left( \frac{ 1 }{ \varkappa_{n} } K^{(n)} \Big(\frac{t }{ \varsigma_{n} } \Big)  - \Big( 1 - t \frac{\alpha}{2} \Big)  \right) \right)_{t \in [0,c]} \Rightarrow \left( \beta_1 - \beta_2 \frac{t}{2} +  B \Big( \alpha \frac t 2 \Big) \right)_{t \in [0,c]}
	\]
	in $D([0,c], \mathbb{R})$ as $n \to \infty$, and 
	where $\big( B(t) \big)_{t \ge 0}$ is a Brownian motion.
\end{thm}

\begin{rem}[On the behavior of \eqref{eq:cond-Thm}]
	Note that the Cauchy--Schwarz inequality implies that
	\[
	(\sigma_n^{(1)})^2 = \left( \sum_{i = 1}^{\kappa_{n}} z_i^{(n)} \right)^2 \le \kappa_{n} \sum_{i = 1}^{\kappa_{n}} \left(z_i^{(n)} \right)^2 = \kappa_{n} \sigma_n^{(2)}.
	\]
	Hence, it is impossible for $(\sigma_n^{(1)})^2 / \kappa_{n} \sigma_n^{(2)}$ (and hence $(\sigma_n^{(1)})^2 / \varkappa_{n} \varsigma_{n}$) to diverge to infinity, and there is always at least a convergent sub-sequence.
	Furthermore, when $\alpha$ as in the statement of Theorem \ref{thm:main_result} exists, it needs to be an element in $[0,1]$.
	However, \eqref{eq:cond-Thm2} requires a stronger assumption on the convergence of $(\sigma_n^{(1)})^2 / \varkappa_{n} \varsigma_{n}$ towards $\alpha$ and the speed of convergence being of order $O(1/\sqrt{\varkappa_{n}})$.
	Besides, since $\kappa_{n} \to \infty$, then in order for $\alpha$ to be strictly positive, at least one of the sequences $\sigma_n^{(1)}$ or $\sigma_n^{(2)}$ needs to have an extreme behavior, i.e.\ $\sigma_n^{(1)} \to \infty$ or $\sigma_n^{(2)} \to 0$.
	
	The utility of scaling by $\varsigma_{n}$ and $\varkappa_{n}$ instead of $\sigma_n^{(2)}$ and $\kappa_{n}$ may not be immediately clear. 
	However, as shown in the examples explored in Section \ref{sec:examples}, for certain models it is more natural or informative to use the equivalent $\varsigma_{n}$ and $\varkappa_{n}$ rather than the actual sequences $\sigma_n^{(2)}$ and $\kappa_{n}$, whose expressions could be more complicated.
\end{rem}

\begin{rem}[On the FCLT] \label{rmk:FCLT}
	The particular version of the martingale Functional Central Limit Theorem (FCLT) that we use in the sequel is the following:
	for each $n \in \mathbb{N}$, let $M_n$ be a square-integrable local martingale in $D([0,T], \mathbb{R})$ satisfying $M_n(0) = 0$.
	Let us assume that the expected maximum jumps of $\langle M_n \rangle$ and $M_n^2$ are asymptotically negligible, i.e.\ for every $T > 0$ we have
	\begin{align*}
		\lim_n \mathbb{E}[J(\langle M_n \rangle , T)] = 0 \text{ and }
		\lim_n \mathbb{E}[J(M_n , T)^2] = 0,
	\end{align*}
	where $J(\phi, T)$ is the absolute value of the maximum jump in $\phi \in D([0,T], \mathbb{R})$ over the interval $[0, T]$, i.e.,
	\[
	J(\phi, T) = \sup \{ |\phi(t) - \phi(t-)| : 0 < t \le T \}.
	\]
	In addition, there exists a $c > 0$ such that for every $t \ge 0$ we get $\langle M_n \rangle_t \Rightarrow c \, t $ as $n \to \infty$, for every $t \in [0,T]$.
	Then, $M_n$ converges toward $(B(c \, t))_{t \in  [0,T]}$ when $n \to \infty$, where $B$ is a Brownian motion.
	
	The statement and proof of this classical result can be found in \cite[Thm.\ 2.1]{Whitt2007}, cf. \ \cite[Thm.\ 1.4 p.\ 339]{EthierKurtz1986}.
	As we commented before, we are able to strengthen the metric of the convergence to the uniform metric because the limit process is continuous.
\end{rem}

\begin{proof}[Proof of Theorem \ref{thm:main_result}]
	According to Proposition \ref{prop:martingale_problem}, scaling the time by $1/\varsigma_{n}$ and the space by $\varkappa_{n}$, the process $M^{(n)} = \left( M^{(n)} \left( t \right) \right)_{t \ge 0}$, where for every $t \ge 0$:
	\[
	M^{(n)} \left( t \right) = \frac{1 }{\varkappa_{n} } K^{(n)} \left( \frac{t}{ \varsigma_{n}} \right) - \frac{\kappa_{n}}{\varkappa_{n}} + \frac{t}{2} \frac{ (\sigma_n^{(1)})^2}{ \varkappa_{n} \varsigma_{n} } - \frac{1}{2 \varkappa_{n} \varsigma_{n} } \int_0^t \sum_{i \ge 1} \Big( X_i^{(n)} \Big)^2 \left( \frac{s}{  \varsigma_{n} } \right) \mathrm{d}s,
	\]
	is a local $(\mathcal{G}_t)$-martingale with predictable quadratic variation
	\[
	\langle M^{(n)} \rangle_t = \frac{t}{2} \frac{ (\sigma_n^{(1)})^2}{ (\varkappa_{n})^2 \varsigma_{n}} - \frac{1}{2 (\varkappa_{n})^2 \varsigma_{n}} \int_0^t \sum_{i \ge 1} \Big( X_i^{(n)} \Big)^2 \left( \frac{s}{  \varsigma_{n} } \right) \mathrm{d}s,
	\]
	where $\mathcal{G}_t = \mathcal{F}_{t / \varsigma_{n}}$, for every $t \ge 0$.
	Notice that
	$$
	\frac{\kappa_{n}}{\varkappa_{n}}  \xrightarrow[n \to \infty]{} 1 \text{ and } \frac{ (\sigma_n^{(1)})^2}{ \varkappa_{n} \varsigma_{n} } \xrightarrow[n \to \infty]{} \alpha.
	$$
	Moreover, we expect the integral term to be negligible, because of Lemma \ref{lemma:bound_2nd-moment_MC}.
	Thus, the limit of $\frac{1 }{\varkappa_{n} } K^{(n)} \left( \frac{t}{ \varsigma_{n}} \right)$ must be $1  - t \frac{\alpha}{2}$.
	Define
	\begin{align*}
		\Delta^{(n)}(t) &:= \frac{1 }{\varkappa_{n} } K^{(n)} \left( \frac{t}{ \varsigma_{n}} \right) -  \left( 1  - t \frac{\alpha}{2} \right).
	\end{align*}
	Then,
	$$
	\Delta^{(n)}(t) = M^{(n)}(t)  + \frac{\kappa_{n}}{\varkappa_{n}} - 1 + \Theta_1^{(n)}(t) + \Theta_2^{(n)}(t),
	$$
	where
	\begin{align*}
		\Theta_1^{(n)}(t) =  \frac{t}{2} \left( \alpha - \frac{ (\sigma_n^{(1)})^2}{ \varkappa_{n} \varsigma_{n}}  \right) \; \text{ and } \;
		\Theta_2^{(n)}(t) = \frac{1}{2 \varkappa_{n} \varsigma_{n} } \int_0^t \sum_{i \ge 1} \Big( X_i^{(n)} \Big)^2 \left( \frac{s}{  \varsigma_{n} } \right) \mathrm{d}s.
	\end{align*}
	Now, $\Theta_1^{(n)}$ converges to zero uniformly on $[0,c]$ because of \eqref{eq:cond-Thm}.
	Furthermore, $\Theta_2^{(n)}$ is non-negative and non-decreasing in $t$, so it suffices to prove the convergence to zero when $t = c$.
	Note that, for $n$ large enough such that $(c / \varsigma_{n}) \cdot \sigma_n^{(2)} < 1$ and we can use the Fubini--Tonelli Theorem, together with Lemma \ref{lemma:bound_2nd-moment_MC}, to obtain
	\begin{equation}\label{eq:control_Theta2}
		\mathbb{E} \left[ \Theta_2^{(n)}(c) \right] \le \frac{1}{2 \varkappa_{n} \varsigma_{n}} \int_0^{c} \frac{ \sigma_n^{(2)} \mathrm{d}s}{ 1 - (s/\varsigma_{n}) \cdot \sigma_n^{(2)}} = - \frac{1}{2 \varkappa_{n}} \ln \left(1 - c \frac{\sigma_n^{(2)} }{ \varsigma_{n} } \right) \xrightarrow[n \to \infty]{} 0,
	\end{equation}
	which implies the convergence in probability uniformly on $[0,c]$ of $\Theta_2$ to zero, as $n \to \infty$.
	Finally, an analogous argument implies that the predictable quadratic variation $\langle M^{(n)} \rangle$ also converges to zero, and then $M^{(n)}$ converges in probability to zero uniformly on $[0,c]$, by using the Doob's maximal inequality or the Burkholder--Davis--Gundy inequality.
	
	Now, to prove the second part of the theorem, note that
	\[
	\sqrt{\varkappa_{n}} \Delta^{(n)}(t) = \sqrt{\varkappa_{n}} M^{(n)}(t)  + \sqrt{\varkappa_{n}} \left( \frac{\kappa_{n}}{\varkappa_{n}} - 1 \right) + \sqrt{\varkappa_{n}} \Theta_1^{(n)}(t) + \sqrt{\varkappa_{n}} \Theta_2^{(n)}(t).
	\]
	Because of \eqref{eq:cond-Thm2}, we have that $t \mapsto  \sqrt{\varkappa_{n}} \left( \frac{\kappa_{n}}{\varkappa_{n}} - 1 \right) +  \sqrt{\varkappa_{n}} \Theta_1^{(n)}(t)$, which is a deterministic process, converges towards $t \mapsto \beta_1 - \beta_2 t /2$ uniformly on $[0,c]$.
	In addition, due to the control of the expectation of $\Theta_2^{(n)}(t)$ established in \eqref{eq:control_Theta2}, we also get that  $\sqrt{\varkappa_{n}} \Theta_2^{(n)}(c) \to 0$.
	
	Note that the jumps of $\sqrt{\varkappa_{n}} M^{(n)}$ are given by the jumps of $\sqrt{\varkappa_{n}} K^{(n)}$, which are all of size  $1/\sqrt{\varkappa_{n}}$ and so they converge to zero.
	Finally, we have for every $t \in [0,c]$:
	\[
	\left\langle \sqrt{\varkappa_{n}} M^{(n)} \right\rangle_t = \frac{t}{2} \frac{ (\sigma_n^{(1)})^2}{ \varkappa_{n} \varsigma_{n}} - \frac{1}{2 \varkappa_{n} \varsigma_{n}} \int_0^t \sum_{i \ge 1} \big( X_i^{(n)} \big)^2 \left( \frac{s}{  \varsigma_{n} } \right) \mathrm{d}s \Rightarrow \frac{t}{2} \alpha, \; \text{ as } n \to \infty
	\]
	and by using the martingale functional central limit theorem (see Remark \ref{rmk:FCLT}) we have that
	\[
	\left( \sqrt{\varkappa_{n}} M^{(n)}(t) \right)_{t \in [0,c]} \Rightarrow \left( B \left(  t \frac{\alpha}{2} \right) \right)_{t \in [0,c]} \text{ as } n \to \infty 
	\]
	concluding the proof.
\end{proof}

In next section we explore the consequence of Theorem \ref{thm:main_result} for some specific examples of random graph processes.

\section{Examples}\label{sec:examples}

Theorem \ref{thm:main_result} is rather general, and it can be applied to several random graph processes.
To illustrate this, we now consider, as examples, the Erd\H{o}s--Rényi random graph (see also Theorem \ref{thm:ER} in the introduction), a non-homogeneous version of this model and a dynamical Norros--Reittu random graph process.

\subsection{Erd\H{o}s--Rényi random graph process}

Let us take $\boldsymbol{z}^{(n)}$ to be the vector with $n$ non-null components all equal to one.
Then, $\sigma_n^{(1)} = \sigma_n^{(2)} = \kappa_n = n$.
Let us denote by $K_{\mathrm{ER}}^{(n)}$ the process recording the number of connected components in $\mathrm{ER}^{(n)}$.
Then, Theorem \ref{thm:ER} in the introduction is immediate from Theorem \ref{thm:main_result}, where $\varsigma_{n} = \varkappa_{n} = n$, $\alpha  = 1$ and $\beta_1 = \beta_2 = 0$.

Our method allows us to study more general models.
One way to add some non-homogeneity to the Erd\H{o}s--Rényi random graph model would be to add a unique vertex of mass $\vartheta > 0$, so that $\boldsymbol{z}^{(n)}$ now consists of $n$ components equal to one, one component equal to $\vartheta$ (located at the beginning or at the end of $\boldsymbol{z}^{(n)}$ depending on whether $\vartheta \ge 1$ or $\vartheta < 1$), and infinitely many zeros after.
It is easy to see that this perturbed model has the same scaling limit as the Erd\H{o}s--Rényi graph process.
In fact, we could even take $\vartheta_{n}$ depending on $n$ and such that $\vartheta_{n}/ \sqrt{n} \to 0$, obtaining the same scaling limit again.

Let us then consider a more general model, where the vector of initial masses $\boldsymbol{z}^{(n)}$ satisfies
\[
\boldsymbol{z}^{(n)} = \mathrm{ord} \big( (\underbrace{1, 1, \dots, 1}_{n \text{ times}}, \vartheta_1, \vartheta_2, \dots, \vartheta_{m_n}, 0, 0, \dots) \big),
\]
and where $\mathrm{ord} : l^2 \to l^2_\searrow$ is the natural projection of $l^2$ onto $l^2_\searrow$.
In words, $\boldsymbol{z}^{(n)}$ consists of $n$ elements equal to $1$, $m_n$ elements respectively equal to $\vartheta_i$ for $i \in \{1, 2, \dots, m_n\}$, and infinitely many zeros after.
Let us denote by $K_{\mathrm{gER}}^{(n)}$ the process recording the number of connected components in this random graph model.

\begin{cor}[Generalized Erd\H{o}s--Rényi random graph]\label{corol:GER}
	
	Assume there exist two non-negative constants $\beta_1$ and $\beta_2$ such that
	\begin{equation}\label{eq:assumpt_Thm}
		\frac{m_n}{\sqrt{n}} \xrightarrow[n \to \infty]{} \beta_1, \; \frac{2}{\sqrt{n}} \sum_{i = 1}^{m_n} \vartheta_i \xrightarrow[n \to \infty]{} \beta_2 \; \text{ and }  \; \frac{1}{n} \sum_{i = 1}^{m_n} \vartheta_i^2 \xrightarrow[n \to \infty]{} 0.
	\end{equation}
	Then, for every $[0,c]$, with $c < 1$ we have that as $n \to \infty$ the sequence of processes
	\[
	\left( \frac{ 1 }{ n } K_{\mathrm{gER}}^{(n)} \Big( \frac{t }{n} \Big) \right)_{t \in [0,c]} \Rightarrow \left( 1 - \frac{t}{2} \right)_{t \in [0,c]} \text{ in } D([0,c], \mathbb{R}) \text{ as } n \to \infty.
	\]
	In addition, the sequence of processes
	\[
	\left( \sqrt{n} \left( \frac{1}{ n }  K_{\mathrm{gER}}^{(n)} \Big( \frac{ t }{ n} \Big)  - \Big( 1 - \frac{t}{2} \Big)  \right) \right)_{t \in [0,c]} \Rightarrow \left( \beta_1 - \beta_2 t + B \Big( \frac t 2 \Big) \right)_{t \in [0,c]}
	\]
	in $D([0,c], \mathbb{R})$ as $n \to \infty$,
	where $B$ is a Brownian motion.
\end{cor}

\begin{rem}
	Notice that \eqref{eq:assumpt_Thm} can be always satisfied by taking, for instance,
	\[
	m_n = \lfloor \beta_1 \sqrt{n} \rfloor \text{ and } \vartheta_i = \frac{\beta_2}{\beta_1}, \text{ for every } i \ge 1. 
	\]
\end{rem}

\begin{proof}
	
	Notice that
	\[
	\frac{ \sigma^{(2)}(\boldsymbol{z}^{(n)}) }{ n } = 1 + \frac{1}{n} \sum_{i = 1}^{m_n} \vartheta_i^2 \xrightarrow[n \to \infty]{} 1 \;\; \text{ and } \;\; \frac{ \kappa(\boldsymbol{z}^{(n)}) }{ n } = 1 + \frac{m_n}{n} \xrightarrow[n \to \infty]{} 1.
	\]
	In addition,
	\[
	\frac{\left( \sigma^{(1)}(\boldsymbol{z}^{(n)}) \right)^2}{n^2} = \frac{1}{n^2} \left( n + \sum_{i = 1}^{m_n} \vartheta_i \right)^2 \xrightarrow[n \to \infty]{} 1.
	\]
	Hence, taking $\varsigma_{n} = \varkappa_{n} = n$ in Theorem \ref{thm:main_result} we immediately get the first part for the result.
	To prove the second part, it suffices to notice that
	\begin{align*}
		\sqrt{n} \left( \frac{n + m_n}{n} - 1 \right) &= \frac{m_n}{\sqrt{n}} \xrightarrow[n \to \infty]{} \beta_1, \text{ and } \\
		\sqrt{n} \left( \frac{1}{n^2} \left( n + \sum_{i = 1}^{m_n} \vartheta_i \right)^2 - 1 \right) &= \frac{ 2 }{\sqrt{n}} \sum_{i = 1}^{m_n} \vartheta_i + \frac{1}{\sqrt{n}} \left( \frac{1}{\sqrt{n}} \sum_{i = 1}^{m_n} \vartheta_i \right)^2 \xrightarrow[n \to \infty]{} \beta_2.
	\end{align*}
\end{proof}

\subsection{A Norros--Reittu random graph process}\label{sec:NR_graph}

Here we consider a version of the multiplicative random graph process, where the weight of each vertex is given by the quantile of a given cumulative distribution function.
Notice this is not the original discrete-time multi-graph model introduced by Norros and Reittu \cite{Norros_Reittu_2006}, which also allows for immigration of vertices and deletion of edges, cf. \cite[\S\ 6.8.2]{Remco2017}.

Let $F$ be a cumulative distribution function of a non-negative random variable $W$, with finite first and second moments.
Let us take
\[
w_i^{(n)} := F^{-1} \left( 1 - \frac i n \right), \text{ for } i \in [n],
\]
where $F^{-1}$ is the generalized inverse of $F$, i.e.
\[
F^{-1}(x) := \inf\{s \in \mathbb{R}_+ : F(s) \ge x \},
\]
and where we set $F^{-1}(0) = 0$.
Let also define
\[
z_i^{(n)} := \frac{w_i^{(n)}}{ \sqrt{l_n}}, \text{ where } l_n = \sum_{i \ge 1} w_i^{(n)}.
\]
Thus, it is clear that $\kappa(\boldsymbol{z}^{(n)}) = n$ and
\begin{align*}
	\frac{\sigma_n^{(1)} }{\sqrt{n}} &= \frac{\sigma^{(1)}(\boldsymbol{z}^{(n)})}{\sqrt{n}} =  \left( \frac{1}{ n} \sum_{i = 1}^n w_i^{(n)} \right)^{1/2} \xrightarrow[n \to \infty]{}  \sqrt{\mathbb{E}[W]} \;\; \text{ and } \\
	\sigma_n^{(2)} &= \sigma^{(2)} (\boldsymbol{z}^{(n)}) = \frac{1}{l_n}\sum_{i = 1}^n \left(w_i^{(n)} \right)^2 \xrightarrow[n \to \infty]{} \frac{\mathbb{E}[W^2]}{\mathbb{E}[W]}.
\end{align*}
The previous two statements are a consequence of the next elementary lemma whose proof, for the sake of completeness, is provided in Appendix \ref{sec:appendix}.

\begin{lem}[Riemann sum of improper integrals]\label{lemma:riemann_sum}
	Let $\phi : (0,1] \to \mathbb{R}_+$ be a non-increasing function such that $\int_0^1 \phi(x) \mathrm{d}x$ is finite, then
	\[
	\lim_n \frac{1}{n} \sum_{k = 1}^n \phi \left( \frac{k}{n} \right) = \int_0^1 \phi(x) \mathrm{d}x.
	\]
	Furthermore, if $\int_0^1 \phi(x)^2 \mathrm{d}x$ is also finite we get
	\[
	\lim_n \sqrt{n} \left(  \int_0^1 \phi(x) \mathrm{d}x - \frac{1}{n} \sum_{k = 1}^n \phi \left( \frac{k}{n} \right)  \right) = 0.
	\]
\end{lem}
Indeed, it suffices to take $\phi$ in Lemma \ref{lemma:riemann_sum} as $x \mapsto F^{-1}(1 - x)$ and $x \mapsto  F^{-1}(1 - x)^2$.
Then, for a uniform on $(0, 1)$ random variable denoted $U$ we have that $F^{-1}(1-U)$ has distribution function $F$.
Hence, we obtain
\[
\int_0^1 F^{-1}(1 - x)^k \mathrm{d}x = \mathbb{E} \left[ \Big( F^{-1}(1 - U) \Big)^k \right] =  \mathbb{E}[W^k].
\]

Let us denote $K_{\mathrm{NR}}^{(n)}$ the process recording the number of connected components in this random graph model starting with $\boldsymbol{z}^{(n)}$ as the initial vector of vertex sizes.

\begin{cor}[Norros--Reittu random graph process]\label{cor:NR}
	For every $[0,c] \subset \left[ 0, \frac{\mathbb{E}[W]}{\mathbb{E}[W^2]} \right)$ we have that the sequence of processes
	\[
	\left( \frac{ 1 }{ n } K_{\mathrm{NR}}^{(n)} \big( t \big) \right)_{t \in [0,c]} \Rightarrow \left( 1 - \frac{t}{2} \mathbb{E}[W] \right)_{t \in [0,c]} \text{ in } D([0,c], \mathbb{R}) \text{ as } n \to \infty.
	\]
	In addition, the sequence of processes of fluctuations
	\[
	\left( \sqrt{n} \left( \frac{1}{ n }  K_{\mathrm{NR}}^{(n)} \big( t \big)  - \Big( 1 - \frac{t}{2} \mathbb{E}[W]  \Big)  \right) \right)_{t \in [0,c]} \Rightarrow \left( B \Big( \frac{t \mathbb{E}[W]}{2} \Big) \right)_{t \in [0,c]}
	\]
	in $D([0,c], \mathbb{R})$  as $n \to \infty$, 
	where $B$ is a Brownian motion.
\end{cor}

\begin{proof}
	Notice that
	\(
	\sigma_n^{(2)} \xrightarrow[n \to \infty]{} \varsigma = {\mathbb{E}[W^2]}/{\mathbb{E}[W]}.
	\)
	Moreover,
	\[
	\frac{(\sigma_n^{(1)})^2}{ \varsigma \, n } = \frac{1}{\varsigma \, n} \sum_{i = 1}^n w_i^{(n)} \xrightarrow[n \to \infty]{} \frac{\mathbb{E}[W]}{\varsigma} = \frac{(\mathbb{E}[W])^2}{\mathbb{E}[W^2]}.
	\]
	Hence, applying Theorem \ref{thm:main_result} we have that for every $[0,c] \subset [0,1]$ the sequence of processes
	\[
	\left( \frac{ 1 }{ n } K_{\mathrm{NR}}^{(n)} \left( \frac{t}{ \varsigma } \right) \right)_{t \in [0,c]} \Rightarrow \left( 1 - \frac{t}{2} \frac{\mathbb{E}[W]}{\varsigma} \right)_{t \in [0,c]} \text{ in } D([0,c], \mathbb{R}) \text{ as } n \to \infty.
	\]
	In addition,
	\[
	\sqrt{n} \left( \frac{1}{\varsigma n} \sum_{i = 1}^n w_i^{(n)} - \frac{\mathbb{E}[W]}{\varsigma}
	\right) = \frac{1}{\varsigma} \cdot \sqrt{n} \left( \frac{1}{n} \sum_{i = 1}^n w_i^{(n)} - \mathbb{E}[W] \right) \xrightarrow[n \to \infty]{} 0,
	\]
	because of Lemma \ref{lemma:riemann_sum}.
	Thus, applying the second half of Theorem \ref{thm:main_result} we have that when $n \to \infty$, the sequence of processes
	\[
	\left( \sqrt{n} \left( \frac{ 1 }{ n } K^{(n)} \Big(\frac{t }{ \varsigma } \Big)  - \Big( 1 -  \frac{t}{2} \frac{\mathbb{E}[W]}{\varsigma}\Big)  \right) \right)_{t \in [0,c]} \Rightarrow \left( B \left( \frac{\mathbb{E}[W]}{\varsigma} \frac t 2 \right) \right)_{t \in [0,c]},
	\]
	in $D([0,c], \mathbb{R})$ as $n \to \infty$.
	The result, as stated in the Corollary \ref{cor:NR} is obtained by the change of variables $s = t / \varsigma$.
\end{proof}

\section{Conclusions and perspectives}

Our methods are quite simple and are based on the use of the martingale problem associated to the multiplicative coalescent dynamics that encodes the evolution of the weighted size of the connected components of the multiplicative random graph process.
It is natural to think that they could be coupled to those of \cite{enriquez2023} and \cite{CorujoLemaireLimic_2024} to study the bi-dimensional process recording the size of the largest component and the number of connected component in the super-critical regime of the Erd\H{o}s--Rényi random graph process.
Even more, including also the method to study the surplus edges of the random graph process developed in \cite{Corujo_Limic_2023a} and \cite{Corujo_Limic_2023b}, we expect to be able to extend our study to cover the joint distribution of three variables: the size of the giant component, the number of surplus edges and the number of connected components, extending the statics results in \cite[Thm.\ 2.2]{Puhalskii2005} to a dynamical setting.
This work can be seen as one step towards that goal.

\appendix

\section{Proof of Lemma \ref{lemma:riemann_sum}}\label{sec:appendix}

Let $U$ be a uniform random variable on $[0,1]$ and let us define the random variable $U_n^{\uparrow}$ as
\begin{equation*}
	U_n^{\uparrow} = \frac{\lceil n U \rceil}{n},
\end{equation*}
so that $U_n^{\uparrow} \ge U$  and $0 \le \phi (U_n^{\uparrow}) \le \phi(U)$.
Note that
\[
\mathbb{E}[\phi(U)] = \int_0^1 \phi(u) \mathrm{d}u < \infty \text{ and } \mathbb{E}[\phi(U_n^{\uparrow})] = \frac{1}{n} \sum_{k = 1}^n \phi \left( \frac{k}{n} \right).
\]
Thus, by using the Dominated Convergence Theorem we get the convergence of $\mathbb{E}[\phi(U_n^{\uparrow})]$ towards $\mathbb{E}[\phi(U)]$, when $n \to \infty$, proving the first part of the result.

For proving the control on the speed of convergence, let us define 
\[
U_n^{\downarrow} = \frac{\lfloor n U \rfloor}{n}, \text{ so that } 0 \le U_n^{\downarrow} \le U \le U_n^{\uparrow}.
\]
In addition, on $(U \ge 1/n)$ we have $0 \le \phi(U_n^{\uparrow}) \le \phi(U) \le \phi(U_n^{\downarrow})$ (notice that outside of this set $\phi(U_n^{\downarrow})$ could be undefined if $\phi$ is undefined at $0$).
Hence,
\begin{align*}
	\mathbb{E} \left[ \sqrt{n} \left( \phi(U) - \phi(U_n^{\uparrow}) \right) \right] \le \mathbb{E} \left[ \sqrt{n} \left( \phi(U_n^{\downarrow}) - \phi(U_n^{\uparrow})\right) \mathbf{1}_{(U \ge \frac 1 n)} \right] + \mathbb{E} \left[ \sqrt{n} \left( \phi(U) - \phi(U_n^{\uparrow})\right) \mathbf{1}_{(U \le \frac 1 n)} \right].
\end{align*}
Notice that, on the one hand, the first term in the RHS of the previous inequality equals
\begin{align*}
	\frac{1}{\sqrt{n}} \left(  \phi \left( \frac{1}{n} \right) - \phi(1) \right) &\le \frac{1}{\sqrt{n}} \phi \left( \frac{1}{n} \right) = \sqrt{ 2 \, \phi^2 \left( \frac{1}{n} \right) \mathbb{E} \left[ \mathbf{1}_{(\frac{1}{2n} \le U \le \frac 1 n)} \right] } \\
	&\le \sqrt{ 2 \, \mathbb{E}\left[\phi^2\left(U \right) \mathbf{1}_{(\frac{1}{2n} \le U \le \frac 1 n)} \right] }
\end{align*}
where the last inequality is a consequence of the fact that $\phi^2$ is decreasing.
Now, the upper bound in the previous chain of inequalities converges to zero because $\mathbb{E}[\phi^2(U)] < \infty$ and the  Dominated Convergence Theorem.
On the other hand, the second term can be bounded as follows
\[
\mathbb{E} \left[ \sqrt{n} \left( \phi(U) - \phi(U_n^{\uparrow})\right) \mathbf{1}_{(U \le \frac 1 n)} \right] \le \mathbb{E} \left[  \sqrt{n} \cdot \phi(U)  \mathbf{1}_{(U \le \frac 1 n)} \right] \le \sqrt{  \mathbb{E} \left[  \phi^2(U)  \mathbf{1}_{(U \le \frac 1 n)} \right] }
\]
where the second inequality is obtained by applying the Cauchy--Schwarz inequality to the product of the random variables $\phi(U) \mathbf{1}_{(U \le \frac 1 n)}$ and $\sqrt{n} \cdot \mathbf{1}_{(U \le \frac 1 n)}$.
Hence, the RHS term converges to zero because $\mathbb{E}[\phi^2(U)] < \infty$ and the  Dominated Convergence Theorem.

\section*{Acknowledgements}

The author would like to thank Nathanaël Enriquez, Gabriel Faraud, Sophie Lemaire, and Vlada Limic for several fruitful discussions that contributed to this work. 
Special thanks are due to Sophie Lemaire and Vlada Limic for their careful reading of the manuscript and their helpful comments. 
The author is also grateful to two anonymous reviewers for their detailed and constructive reports.


\providecommand{\bysame}{\leavevmode\hbox to3em{\hrulefill}\thinspace}
\providecommand{\MR}{\relax\ifhmode\unskip\space\fi MR }
\providecommand{\MRhref}[2]{%
	\href{http://www.ams.org/mathscinet-getitem?mr=#1}{#2}
}
\providecommand{\href}[2]{#2}

\end{document}